%% file: changcof.tex
\documentclass[a4paper,english]{scrartcl}
\usepackage[utf8x]{inputenc}
\usepackage[T1]{fontenc}
\usepackage{babel}
\usepackage{amsmath,amssymb}
\usepackage[amsmath,thref,thmmarks]{ntheorem}
\usepackage[pdftex]{graphicx}
\usepackage{color}
\usepackage{url}
\usepackage[all]{xy}

\input makros

\begin{document}

\title{Some basic thoughts on the cofinalities of Chang Structures with an application to forcing}
\author{Dominik T. Adolf}
\maketitle

\begin{abstract}
 Consider $(\kappa^{+++},\kappa^{++}) \twoheadrightarrow (\kappa^+,\kappa)$ where $\kappa$ is an uncountable regular cardinal. By a result of Shelah's we have $\cof(X \cap \kappa^{++}) = \kappa$ for almost all $X \subset \kappa^{+++}$ witnessing this. Here we consider the question if there could be a similar result for $X \cap \kappa^{+}$. We will use this discussion to give an interesting example of an pseudo Prikry forcing answering a question of Sinapova.
\end{abstract}

\section{Introduction}

Let $\lambda,\lambda',\kappa,\kappa'$ be regular cardinals with $\lambda > \kappa$ and $\lambda' > \kappa'$. We take $(\lambda,\kappa) \twoheadrightarrow (\lambda',\kappa')$ to mean: for any structure (in a countable language) on $\lambda$ there exists a substructure $X$ with $\card(X) = \lambda'$ but $\card(X \cap \kappa) = \kappa'$. 

Originating in Model Theory this seemingly innocuous property has significant large cardinal strength. Its most basic form (often known simply as \textit{the} Chang conjecture), $(\aleph_2,\aleph_1) \twoheadrightarrow (\aleph_1,\aleph_0)$ in our notation, is equiconsistent with an $\omega_1$-Erd\H{o}s cardinal.

The exact consistency strength of the analogous property $(\aleph_3,\aleph_2) \twoheadrightarrow (\aleph_2,\aleph_1)$ has proved more elusive. Though it is known to be consistent relative to the existence of a huge cardinal by an argument of Kunen's \cite{kunensatideals}.

Generally, we will use the language of stationary sets throughout this paper. Recall, a subset $S$ of $\Pot(\mathcal{X})$ is stationary iff for all $F:\finsubsets{\mathcal{X}} \rightarrow \mathcal{X}$ there is some $X \subseteq \mathcal{X}$ in $S$ that is closed under $F$. $(\lambda,\kappa) \twoheadrightarrow (\lambda',\kappa')$ then means: there is a stationary set of $X \subset \lambda$ with $\card(X) = \lambda'$ but $\card(X \cap \kappa) = \kappa'$. 

Note that this is equivalent to saying that there is stationary set of $X \subset H_\lambda$ with $X \cap \lambda$ satisfying $(\lambda,\kappa) \twoheadrightarrow (\lambda',\kappa')$.

Consider now some $\kappa$ uncountable and regular. Assume $(\kappa^{+++},\kappa^{++}) \twoheadrightarrow (\kappa^+,\kappa)$. (A more optimal consistency proof for this property has recently appeared in\cite{ccfromsubcomp}). The question given $X \subset \kappa^{+++}$ witnessing this and a cardinal $\alpha \leq \kappa^{+++}$ is: what is $\cof(X \cap \alpha)$?

To avoid unnecessary complexity we will make the additional assumption that $\kappa \subset X$. (The known consistency proofs provide structures with this property. This stronger property will be denoted by the subscript $\kappa$ as in $(\kappa^{+++},\kappa^{++}) \twoheadrightarrow_\kappa (\kappa^+,\kappa)$) Furthermore, due to an argument of Shelah's (see \cite{Jech}, p 451, first applied in this context in \cite{fandmdefcounterexampletoch}) we have $\cof(X \cap \kappa^{++}) = \kappa$ for all but non-stationarily many Chang structures $X$. We remain curious about $\cof(X \cap \kappa^+)$. 

In the first part of this paper we will give some rather basic results about possible values for this cofinality using elementary methods. Most pertinent here is that if $\GCH$ holds, then all reasonable values are possible. 

These results also extend to many other Chang properties, say $(\kappa^{+4},\kappa^{+++}) \twoheadrightarrow_\kappa (\kappa^+,\kappa)$, but to keep our notation tidy we will consider only the simplest representative case.

We will use these results together with arguments from \cite{ccfromsubcomp} to prove this:

\begin{thm}
   Let $\kappa < \lambda < \delta$ be three cardinals such that $\kappa$ is $\lambda$-supercompact, $\lambda$ is $+2$-subcompact, and $\delta$ is Woodin or strongly compact. Then there exists a set generic extension of the universe $V\left[G\right]$ containing partially ordered sets $\partord_\alpha$ for any regular cardinal $\alpha \leq \kappa$ with the following properties: $\partord_\alpha$ does not add bounded subsets to $\kappa$, it changes $\cof(\kappa)$ to $\omega$ and $\cof(\kappa^+)$ to $\alpha$, but does not collapse $\kappa^{+++}$.
\end{thm}

 This answers a question of Sinapova \cite[Question 3]{singcardsurvey}, but we do think it an incomplete answer as the forcings we get do collapse many cardinals in between $\kappa^{+++}$ and $\delta$. The argument does point to a possible approach to a complete answer, though. \thref{fragedrei}

\section{Basic thoughts}

For the remainder of this section we shall fix an uncountable regular cardinal $\kappa$. 

\begin{lemma}
	Assume $(\kappa^{+++},\kappa^{++}) \twoheadrightarrow_\kappa (\kappa^+,\kappa)$ and $2^{\kappa^+} \leq \kappa^{++}$. Then for all regular cardinals $\alpha$ less than $\kappa^+$ we have a stationary set $S_\alpha$ witnessing $(\kappa^{+++},\kappa^{++}) \twoheadrightarrow_\kappa (\kappa^+,\kappa)$ such that $\cof(X \cap \kappa^+) = \alpha$ for all $X \in S_\alpha$.
\end{lemma}

\begin{beweis}
	Let $\mathfrak{A}$ be a Skolemized structure on $H_{\kappa^{+++}}$. By assumption, we know that there exists some $X \prec \mathfrak{A}$ with $\kappa \subset X$, $\card(X \cap \kappa^{+++}) = \kappa^+$, and $\card(X \cap \kappa^{++}) = \kappa$. We shall show that for all $\beta \in \kappa^{+}$ the Skolemhull of $X \cup \{\beta\}$ (notated here as $\Sk^{\mathfrak{A}}$), too, is a Chang structure of the required type. By iterating this procedure we can then construct an $\subset$-increasing and continuous sequence of Chang structures $\<X_\gamma:\gamma \leq \kappa\>$ with $X_\gamma \prec \mathfrak{A}$ and $X_\beta \cap \kappa^+ \in X_\gamma$ for all $\beta < \gamma \leq \kappa$. 
	
	$X_\alpha$ then is a substructure of $\mathfrak{A}$ of the required type with $\cof(X_\alpha \cap \kappa^+) = \alpha$. 
	
	Consider then $Y := \Sk^\mathfrak{A}(X \cup \{\beta\})$. Obviously, we have $\kappa \subset Y$ and $\card(Y \cap \kappa^{+++}) = \kappa^+$. It remains to be seen that $\card(Y \cap \kappa^{++}) = \kappa$. Without loss of generality we can assume that $\mathfrak{A}$ satisfies the collection scheme. Therefore we have that any element of $Y \cap \kappa^{++}$ can be written as $f(\beta)$ where $\dom(f) = \kappa^+$,$\ran(f) \subset \kappa^{++}$ and $f \in X$.
	
	By our assumption there are only $\kappa^{++}$ many functions from $\kappa^+$ into $\kappa^{++}$. We then have some bijective $g: \kappa^{++} \rightarrow \fnktsraum{\kappa^+}{\kappa^{++}}$ in $H_{\kappa^{+++}}$. By elementarity there is some such $g \in X$ witnessing that $\card(\fnktsraum{\kappa^+}{\kappa^{++}} \cap X) = \card(X \cap \kappa^{++}) = \kappa$. We conclude then that $\card(Y \cap \kappa^{++}) = \kappa$ as desired.
\end{beweis}

The next lemma shows that the conclusion of the preceding lemma can hold even if the powerset of $\kappa^+$ is large.

\begin{lemma}
	Assume the conclusion of the preceding lemma holds. Then there exists some forcing extension of the universe $V\left[G\right]$ in which the conclusion remains true, but $2^{\kappa^+} > \kappa^{++}$.
\end{lemma}

\begin{beweis}
	Let $\partord$ be the partial order for adding $\kappa^{+++}$ many Cohen reals. It is well known that $\partord$ has the c.c.c. From this it follows that $S_\alpha$ from the conclusion of the lemma will remain stationary in any $\partord$-generic extension of the Universe.
	
	To see this, realize that for any name $\dot{F}$ for a function in $\fnktsraum{\finsubsets{\kappa^{+++}}}{\kappa^{+++}}$ there exists in $V$ some $F:\omega \times \finsubsets{\kappa^{+++}} \rightarrow \kappa^{+++}$ such that $\dot{F}^G(a) \in \{F(n,a)\vert n < \omega\}$ for all finite $a \subset \kappa^{+++}$ and all $G \subset \partord$ generic over $V$.
\end{beweis}

 So, so far we have realized that we do have complete freedom in what cofinalities can be realized in a Chang structure assuming $\GCH$. We have also realized that $\GCH$ is not strictly required for that conclusion. The next lemma will show though that some restrictions on the behaviour of the continuum function do exist.

\begin{lemma}\label{barringcofinalities}
	Assume $(\kappa^{+++},\kappa^{++}) \twoheadrightarrow_\kappa (\kappa^+,\kappa)$, and $2^\kappa \leq \kappa^{++}$ but $2^{\kappa^+} > \kappa^{++}$. Let $\mu$ be such that $\kappa^\mu = \kappa$, then $\cof(X \cap \kappa^+) > \mu$ for all but non-stationarily many $X \subset \kappa^{+++}$ witnessing $(\kappa^{+++},\kappa^{++}) \twoheadrightarrow_\kappa (\kappa^+,\kappa)$.
\end{lemma}

\begin{beweis}
	Assume for a contradiction that some $X \subset H_{\kappa^{+++}}$ witnesses $(\kappa^{+++},\kappa^{++}) \twoheadrightarrow_\kappa (\kappa^+,\kappa)$ but $cof(X \cap \kappa^+) \leq \mu$. We can then assume that $\card(\Pot(\kappa^+) \cap X) = \card(X \cap \kappa^{+++}) = \kappa^+$ but $\card(\Pot(\kappa) \cap X) = \card(\kappa^{++} \cap X) = \kappa$. 
	
	On the other hand, every subset of $\kappa^+$ in $X$ can be identified with a $\mu$-long sequence of bounded subsets of $\kappa^+$ all of which are in $X$. (Say $\<\alpha_\xi:\xi < \mu\>$ were cofinal in $X \cap \kappa^+$ then any $A \in X \cap \Pot(\kappa^+)$ is determined by $\<A \cap \alpha_\xi: \xi < \mu\>$.)
	
	We then have that $\kappa^+ = \kappa^\mu = \kappa$. Contradiction!
\end{beweis}

\section{A pseudo Prikry forcing}

The term Prikry forcing is here used informally to refer to a forcing notion that changes the cofinality of some (inaccessible) cardinal $\kappa$ but does not add bounded subsets of $\kappa$. Using supercompact cardinals it is possible to construct a forcing that will change the cofinality of both $\kappa$ and any number of succesor cardinals. In such a construction however $\kappa$ and its successors must be assigned one uniform cofinality. Generally, such Prikry forcings will also have optimal chain conditions. 

In \cite{singcardsurvey} Sinapova asked if there could be a forcing extension of the universe $V\left[G\right]$ that changes the cofinalities of say $\kappa,\kappa^+,\kappa^{++}$ without adding bounded subsets, and satisfies that $\cof^{V\left[G\right]}(\kappa) \neq \cof^{V\left[G\right]}(\kappa^+)$ but $\kappa^{+++}$ is not collapsed. (Note it is impossible to have a forcing extension such that $\cof^{V\left[G\right]}(\kappa) \neq \cof^{V\left[G\right]}(\kappa^+)$ but $\kappa^{++}$ is not collapsed.)

Here we will show that the answer is yes. We will have a forcing extension $V\left[G\right]$ such that $\cof^{V\left[G\right]}(\kappa) \neq \cof^{V\left[G\right]}(\kappa^+)$ and $\kappa^{+++}$ is preserved (, but larger cardinals may be collapsed). Our forcing will not add bounded subsets of $\kappa$ but we call it a pseudo Prikry forcing because it does not have the structural properties of a ``true" Prikry forcing.

\begin{lemma}
   Let $\kappa < \lambda < \delta$ be cardinals such that $\kappa$ is $\lambda$-super compact, $\lambda$ is $+2$-subcompact, and $\delta$ is Woodin or strongly compact. Assume $\GCH$. There exists a forcing extension $V\left[H\right]$ such that in $V\left[H\right]$ $\kappa$ is still measurable, $\GCH$ is preserved, and $(\kappa^{+++},\kappa^{++}) \twoheadrightarrow_\kappa (\kappa^+,\kappa)$ holds. Furthermore, $\delta$ remains Woodin or strongly compact.
\end{lemma}

\begin{beweis}

By \cite{ccfromsubcomp} there exists some $\rho < \lambda$ such that forcing with $\col(\kappa,\rho^{+}) \times \col(\rho^{++},\lambda)$ gives us $(\kappa^{+++},\kappa^{++}) \twoheadrightarrow_\kappa (\kappa^+,\kappa)$. We claim that more is true. There is in fact a $\rho < \lambda$ such that forcing with $\partord \ast (\col(\kappa,\rho^{+}) \times \col(\rho^{++},\lambda))^{V^{\partord}}$ gives us $(\kappa^{+++},\kappa^{++}) \twoheadrightarrow_\kappa (\kappa^+,\kappa)$ for any $\partord \subseteq V_\kappa$.
	
	This requires only minor alterations to the proof of Theorem 17 from \cite{ccfromsubcomp}. Counterexamples will now be pairs $\<\partord_\rho,\ddot{f}_\rho\>$ where $\partord_\rho \subset V_\kappa$ is a partial order, and $\ddot{f}_\rho$ is a $\partord_\rho$-name for $\dot{f}_\rho$, a $\col(\kappa,\rho^+) \times \col(\rho^{++},\lambda)$-name for a counterexample to $(\kappa^{+++},\kappa^{++}) \twoheadrightarrow_\kappa (\kappa^+,\kappa)$. We can then pick some $\rho$ such that there exists an elementary embedding
	\[ j: (H(\rho^{++});\in,\<\<\mathbb{Q}_\xi,\ddot{g}_\xi\>:\xi < \rho\>) \rightarrow (H(\lambda^{++});\in,\<\partord_\rho,\ddot{f}_\rho\>:\rho < \lambda\>). \]
	The important part is that if $G \subset \partord_\rho$ is generic over $V$ then $j$ extends onto $H(\rho^{++})$. The rest is as before.
	
	Pick then such a good $\rho$. Pick also some $f:\kappa \rightarrow \kappa$, $g:\kappa \rightarrow \kappa$, and a $\lambda$ supercompactness embedding such that $j(f)(\kappa) = \rho$ and $j(g)(\kappa) = \lambda$. Let $\partord_\kappa$ be the Easton support iteration made up of $\<\partord_\alpha,\dot{Q}_\alpha:\alpha < \kappa\>$ where $\forces_{\partord_\alpha} \dot{Q}_\alpha = \col(\check{\alpha},\check{f(\alpha)}^+) \times \col(\check{f(\alpha)}^{++},\check{g(\alpha)})$ for all $\alpha < \kappa$. Let $G_\kappa \subset \partord_\kappa$ be generic over $V$ and let $\partord := (\col(\kappa,\rho^{+}) \times \col(\rho^{++},\lambda))^{V\left[G\right]}$. Let then $G \subset \partord$ be generic over $V\left[G_\kappa\right]$.
	
	Because $\rho$ was good and $\partord_\kappa \subset V_\kappa$ we do know that $(\kappa^{+++},\kappa^{++}) \twoheadrightarrow_\kappa (\kappa^+,\kappa)$ holds in $V\left[G_\kappa\right]\left[G\right]$. But we would also like to know that $\kappa$ remains measurable. The following is standard: $\partord \concat \dot{P}$ is an initial segment of $j(\partord_\lambda \concat \dot{P})$; $p := \Union\ptwimg{j}{G} \in M\left[G_\kappa\right]\left[G\right]$ and acts as master condition, i.e. if $H_{j(\kappa)} \subset j(\partord)$ is generic over $M\left[G_\kappa\right]\left[G\right]$ end-extending $G_\kappa \concat G$ and $H \subset j(\partord)$ is generic over $M\left[H_{j(\kappa)}\right]$ with $p \in H$, then $j$ extends to an elementary embedding from $V\left[G_\kappa\right]\left[G\right]$ into $M\left[H_{j(\kappa)}\right]\left[H\right]$. The last thing to note is that such a $H$ can be found in $V\left[G_\kappa\right]\left[G\right]$ due to the closure of the collapse forcing and the small amount of $M\left[G_\kappa\right]\left[G\right]$-dense sets.

\end{beweis}

This is the first ingredient in our proof. The second and third ingredients will be the results from the previous section and Woodin's Stationary Tower Forcing.

\begin{lemma}
	Let $\kappa$ be a measurable cardinal such that $(\kappa^{+++},\kappa^{++}) \twoheadrightarrow_\kappa (\kappa^+,\kappa)$ holds. Let $\delta > \kappa$ be  a Woodin cardinal or a strongly compact cardinal. Assume $\GCH$. Fix $\alpha \leq \kappa$ a regular cardinal. Then there exists a forcing notion that will not add bounded subsets to $\kappa$, changes the cofinality of $\kappa$ to $\omega$, that of $\kappa^+$ to $\alpha$ but does not collapse $\kappa^{+++}$.
\end{lemma}

\begin{beweis}
	Let $\partord_\delta$ be the full stationary tower on $V_\delta$ (see \cite{Larson}). By $(\kappa^{+++},\kappa^{++}) \twoheadrightarrow_\kappa (\kappa^+,\kappa)$ and the results of the previous section there exists a stationary set $S_\alpha$ on $\kappa^{+++}$ with the property that $\card(X) = \kappa^+$, $\card(X \cap \kappa^{++}) = \kappa$, $\kappa \subset X$ and $\cof(X \cap \kappa^+) = \alpha$ for all $X \in S_\alpha$. Note that $S_\alpha$ is a condition in $\partord_\delta$.
	
	Let $G \subset \partord_\delta$ be a $V$-generic filter with $S_\alpha \in G$. In $V\left[G\right]$ there exists a generic elementary $j:V \rightarrow M$ with the following properties:
	
	\begin{itemize}
	 \item[$(i)$] $M$ is transitive, in fact, $\fnktsraum{\kleiner\delta}{M} \subset M$;
	 \item[$(ii)$] $S \in G \Gdw \ptwimg{j}{\Union S} \in j(S)$ for all conditions $S \in \partord_\delta$.
	\end{itemize}
	
	From $(ii)$ and the fact that $S_\alpha \in G$ we can then conclude the following facts:
	
	\begin{itemize}
	 \item $\crit(j) = (\kappa^+)^V$,
	 \item $j((\kappa^+)^V) = (\kappa^+)^{V\left[G\right]} = (\kappa^{+++})^V$,
	 \item $\cof^{V\left[G\right]}((\kappa^+)^V) = \alpha$.
	\end{itemize}
	
	 By elementarity $\kappa$ is then measurable in $M$ and hence in $V\left[G\right]$. (The forcing does not add bounded subsets to $\kappa$ because $\kappa$ is inaccessible.)
	
	We can then use some normal measure on $\kappa$ in $V\left[G\right]$ to get a generic extension $V\left[G\right]\left[H\right]$ where $\cof^{V\left[G\right]\left[H\right]}(\kappa) = \omega$ and $\kappa^+$ (formerly $\kappa^{+++}$) is not collapsed. This shows that $(\partord_\delta)_{S_\alpha} \ast \dot{Q}$ ($\dot{Q}$ a $\partord_\delta$-name for the Prikry forcing) is as desired.
\end{beweis}

\section{Open Questions}

We finish on a review of some interesting open problems. First, consider \thref{barringcofinalities}; this lemma could be used to eliminate certain cofinalities from appearing in Chang structures, but it is unclear if it is even consistent.

\begin{frage}
  Is it consistent to have for some regular cardinal $\kappa$, $(\kappa^{+++},\kappa^{++}) \twoheadrightarrow_\kappa (\kappa^+,\kappa)$, and $2^\kappa \leq \kappa^{++}$ but $2^{\kappa^+} > \kappa^{++}$?
\end{frage}

Most interesting would be the ability to do away with countable cofinalities. Ideally, we would like to prove that it is possible to always find Chang structures with uncountable cofinality, assuming of course that we can find Chang structures at all. This would yield better lower bounds for consistency strength, e.g. following \cite{seancc} we could conclude that $(\kappa^{+++},\kappa^{++}) \twoheadrightarrow_\kappa (\kappa^+,\kappa)$ implies the existence of an inner model with a repeat point.

\begin{frage}
	Let $\kappa$ be an uncountable, regular cardinal. Assume $(\kappa^{+++},\kappa^{++}) \twoheadrightarrow_\kappa (\kappa^+,\kappa)$. Does there exist a stationary set of $X \subset \kappa^{+++}$ witnessing $(\kappa^{+++},\kappa^{++}) \twoheadrightarrow_\kappa (\kappa^+,\kappa)$ with the property that $cof(X \cap \kappa^+) > \omega$?
\end{frage}

Finally, our use of the Stationary Tower in our forcing construction has two major downsides. Firstly, we do require an additional large cardinal that exists far above the area we are generally interested in. Secondly, the Stationary Tower will collapse many more cardinals than seemingly necessary for this purpose. We would like to know if it is possible to replace it with a smaller tower of ideals or even a single ideal.

\begin{frage}\label{fragedrei}
  Is it consistent to have $\mu$, a measurable cardinal, such that for all regular $\alpha \leq \kappa$ there exists a presaturated ideal $I_\alpha$ on $\Pot_{\mu^+}(\mu^{++})$ such that $\{ X \in \Pot_{\mu^+}(\mu^{++}) \vert \cof(X \cap \mu^+) = \alpha \}$ is $I_\alpha$-positive?
\end{frage}

\bibliographystyle{alpha}
\bibliography{bibliographie}

\end{document}

%% file: makros.tex
\theoremstyle{plain}
\theoremheaderfont{\normalfont\bfseries}
\theorembodyfont{\itshape}
\theoremseparator{:}
\theoremnumbering{arabic}
\theoremsymbol{}

\newtheorem{thm}{Theorem}[section]

\theoremstyle{plain}
\theoremheaderfont{\normalfont\bfseries}
\theorembodyfont{\itshape}
\theoremseparator{:}
\theoremnumbering{arabic}
\theoremsymbol{}

\theoremstyle{plain}
\theoremheaderfont{\normalfont\bfseries}
\theorembodyfont{\itshape}
\theoremseparator{:}
\theoremnumbering{arabic}
\theoremsymbol{}

\theoremstyle{plain}
\theoremheaderfont{\normalfont\bfseries}
\theorembodyfont{\itshape}
\theoremseparator{:}
\theoremnumbering{arabic}
\theoremsymbol{}

\newtheorem{lemma}[thm]{Lemma}

\theoremstyle{plain}
\theoremheaderfont{\normalfont\bfseries}
\theorembodyfont{\itshape}
\theoremseparator{:}
\theoremnumbering{arabic}
\theoremsymbol{}

\theoremstyle{plain}
\theoremheaderfont{\normalfont\bfseries}
\theorembodyfont{\itshape}
\theoremseparator{:}
\theoremnumbering{Alph}
\theoremsymbol{}

\theoremstyle{plain}
\theoremheaderfont{\normalfont\bfseries}
\theorembodyfont{\itshape}
\theoremseparator{:}
\theoremnumbering{Alph}
\theoremsymbol{}

\theoremstyle{plain}
\theoremheaderfont{\bfseries}
\theorembodyfont{\normalfont}
\theoremseparator{:}
\theoremnumbering{arabic}
\theoremsymbol{}

\theoremstyle{plain}
\theoremheaderfont{\bfseries}
\theorembodyfont{\normalfont}
\theoremseparator{:}
\theoremnumbering{arabic}
\theoremsymbol{}

\theoremstyle{plain}
\theoremheaderfont{\bfseries}
\theorembodyfont{\normalfont}
\theoremseparator{:}
\theoremnumbering{arabic}
\theoremsymbol{}

\theoremstyle{nonumberplain}
\theoremheaderfont{\bfseries}
\theorembodyfont{\normalfont}
\theoremseparator{:}
\theoremnumbering{arabic}
\theoremsymbol{}

\theoremstyle{nonumberplain}
\theoremheaderfont{\bfseries}
\theorembodyfont{\normalfont}
\theoremseparator{:}
\theoremnumbering{arabic}
\theoremsymbol{}

\theoremstyle{plain}
\theoremheaderfont{\normalfont\scshape}
\theorembodyfont{\itshape}
\theoremseparator{:}
\theoremnumbering{arabic}
\theoremsymbol{}

\newcounter{QuesCounter}
\theoremstyle{plain}
\theoremheaderfont{\normalfont\scshape}
\theorembodyfont{\itshape}
\theoremseparator{:}
\theoremnumbering{arabic}
\theoremsymbol{}

\newtheorem{frage}[QuesCounter]{Question}

\theoremstyle{nonumberplain}
\theoremheaderfont{\normalfont\scshape}
\theorembodyfont{\itshape}
\theoremseparator{:}
\theoremnumbering{arabic}
\theoremsymbol{}

\theoremstyle{nonumberplain}
\theoremseparator{:}
\theoremheaderfont{\scshape}
\theorembodyfont{\upshape}
\theoremsymbol{{\Large\ensuremath{_\dashv}}}

\newtheorem{beweis}{Proof}

\theoremstyle{nonumberplain}
\theoremseparator{:}
\theoremheaderfont{\scshape}
\theorembodyfont{\upshape}
\theoremsymbol{\ensuremath{\square}}

\newcommand{\<}{\langle}
\renewcommand{\>}{\rangle}
\renewcommand{\phi}{\varphi}
\newcommand{\Gdw}{\Leftrightarrow}

\newcommand{\op}[1]{\operatorname{#1}}
\newcommand{\Pot}{\op{\mathcal{P}}}

\newcommand{\GCH}{\op{GCH}}

\newcommand{\Union}{\bigcup}

\newcommand{\concat}{{}^\smallfrown}

\newcommand{\dom}{\op{dom}}

\newcommand{\crit}{\op{crit}}
\newcommand{\ran}{\op{ran}}

\newcommand{\col}{\op{Col}}
\newcommand{\card}{\op{Card}}
\newcommand{\cof}{\op{cof}}

\newcommand{\Sk}{\op{Sk}}

\newcommand{\fnktsraum}[2]{{}^{#1}{}_{#2}}
\newcommand{\partord}{\mathbb{P}}
\newcommand{\ptwimg}[2]{{#1}"\left[{#2}\right]}
\newcommand{\kleiner}{\mathord{<}}

\newcommand{\forces}{\Vdash}

\newcommand{\finsubsets}[1]{{\left[#1\right]}^{\kleiner\omega}}